\documentclass[12pt]{amsart}
\usepackage{latexsym}
\usepackage{amsthm}
\usepackage{amssymb}
\usepackage{amsmath}
\usepackage{bm}

\newtheorem{thm}{Theorem}[section]
\newtheorem{lem}[thm]{Lemma}
\newtheorem{con}[thm]{Conjecture}
\newtheorem{cor}[thm]{Corollary}
\newtheorem{prop}[thm]{Proposition}

\theoremstyle{remark}
    \newtheorem*{rem}{Remark}

\theoremstyle{definition}
    \newtheorem*{defn}{Definition}

\numberwithin{equation}{section}

\def\Z {{\mathbb Z}}
\def\N {{\mathbb N}}

\def\H {{\mathcal H}}

\def\C {{\mathcal C}}
\def\l {\langle}
\def\r {\rangle}
\def\ul {\mathbf}

\def\min {{\mathrm{min}}} 

\def\trace {\operatorname{tr}}
\def\LL {\mathcal L}

\makeglossary

\begin{document}
\title[the Dipper-James Conjecture]{Centres of Hecke algebras: \\ the Dipper-James conjecture}
\author{Andrew R. Francis}
\address{School of Computing and Mathematics, University of Western Sydney, NSW 1797, Australia}
\email{a.francis@uws.edu.au}
\author{John J. Graham}
\address{School of Mathematics, University of Sydney, NSW 2006, Australia}
\email{john@xq.se}
\date{15th May 2006}
\begin{abstract}
In this paper we prove the Dipper-James conjecture that the centre of the Iwahori-Hecke algebra
of type $A$ is the set of symmetric polynomials in the Jucys-Murphy operators.
\end{abstract}
\thanks{Thanks to the Courthouse and the Carlisle Castle Hotels in Newtown for their work-friendly atmospheres.}
\maketitle

\section{Introduction}
The symmetric group $W=\mathfrak S_n$ is a Coxeter group generated by the set
$S$ of simple reflections $s_i:=(i,i+1)$ $(1\le i<n)$.
If $w\in\mathfrak S_n$, a product $s_{i_1}s_{i_2}\cdots s_{i_k}=w$ is called reduced when $k$ is minimal.
In this case $k$ is called the length $\ell(w)$ of $w$. \glossary{\ell{w}}

Let $R$ be a commutative ring with $1$ and $q\in R$ be invertible.
The Hecke algebra $\H=\H_n(R,q)$ of the symmetric group is the associative $R$-algebra with basis $T_w$ $(w\in W)$,
and relations induced by the following:
\begin{enumerate}
\item
If $\ell(x)+\ell(y)=\ell(xy)$, then $T_xT_y=T_{xy}$.
\item
If $s\in S$, then $(T_s+1)(T_s-q)=0$.
\end{enumerate}
The algebra is generated by the family $T_i:=T_{s_i}$ $(s_i\in S)$.
The following elements $L_i$ of $\H$ are called Jucys-Murphy operators \cite{Jucys,Murphy}:
$$
L_i:=\sum_{1\le j<i} q^{j-i} T_{(j,i)}\quad (1\le i\le n).
$$
However we find $\mathcal L_i:=(q-1)L_i + 1$ easier to work with because of the recurrence:
$$
\mathcal L_1 = 1\quad\text{and}\quad q\mathcal L_{i+1}=T_i\LL_iT_i.
$$
The symmetric polynomials in $L_1,\dots,L_n$ are central in $\H$
because each generator $T_j$ commutes with $\mathcal L_j+\mathcal L_{j+1}$,
$\mathcal L_j\mathcal L_{j+1}$ and $\mathcal L_i$ $(i\ne j, j+1)$.

Dipper and James have conjectured \cite[Thm 2.41]{DJ}:
\begin{con}[Dipper-James]\label{thm:djconjst}
The centre of the Hecke algebra of the symmetric group
is the set of symmetric polynomials in the Jucys-Murphy operators.
\end{con}
In this paper, we prove this conjecture with $R$ and $q$ as above.
The symmetric group case ($q=1$) over a field was established \cite{Murphy} by Murphy.
In \cite{DJ}, Dipper and James generalise Murphy's result to the Hecke algebra,
but the proof has a gap in the non-semisimple case.

Let $\widehat \H=\widehat \H_n(R,q)$ denote the affine Hecke algebra \cite{L89} associated with the general linear group
over a non-archimedian field.
Fix weight lattice $P=\Z\epsilon_1+\dots+\Z\epsilon_n$ and the ``geometric choice'' for the positive
root system (simple roots $\alpha_i:=\epsilon_{i+1}-\epsilon_i$).
Then $\widehat \H$ has a (Bernstein) basis $\{X_\lambda T_w\}$ indexed by $\lambda\in P$ and $w\in W$.
The span $R[T]=\langle X_\lambda\mid \lambda\in P\rangle$ is
the ring of Laurent polynomials in $X_{\epsilon_i}$ over $R$.
The centre $\widehat Z$ of $\widehat \H$ has been characterised by Bernstein and Zelevinski \cite{L89,Ch}
as the set $R[T]^W$ of symmetric Laurent polynomials.

The action of the large abelian subalgebra $R[T]$ on ``standard'' $\widehat \H$-modules admits a combinatorial description in terms
of Young tableaux. The Specht modules of the Hecke algebra inherit an $\widehat \H$-action via the well known \cite{Ariki,RR}
surjective algebra homomorphism $\psi:\widehat \H\to \H$
which maps $T_w\mapsto T_w$ $(w\in W)$ and $X_{\epsilon_i}\mapsto \mathcal L_i$ $(1\le i\le n)$.
It is immediate that the image of the centre of $\widehat \H$ is contained in the centre of $\H$.
If $q-1$ is invertible, Conjecture~\ref{thm:djconjst} implies that these subsets of $\H$ coincide.

The Hecke algebra is a symmetric algebra with respect to the trace
$$
\trace\left(\sum_{w\in W} a_w T_w\right)=a_1.
$$
The associated bilinear form
$$
\langle T_x, T_y\rangle=\trace(T_yT_x)=
\begin{cases}
q^{\ell(x)}&\text{ if $xy=1$,}\\
0&\text{ otherwise,}
\end{cases}
$$
induces an isomorphism $x\mapsto \langle -,x\rangle$ of $R$-modules between $\H$ and its dual.
The centre $Z$ of $\H$ maps to the space of trace functions,
that is, linear functions $\chi:\H\to R$ such that $\chi(ab)=\chi(ba)$ if $a,b\in \H$.
Geck and Rouquier \cite{GR97,Fmb} have constructed a basis
$$
\{f_{\mathcal C}\mid \mathcal C\hbox{ is a conjugacy class of }W\}
$$
for the space of trace functions (and hence the centre of $\H$).
\begin{thm}\label{thm:basis}\cite[Thm 8.2.3]{GP}
For each conjugacy class $\mathcal C$ of $W$,
there exists a unique trace function $f_{\mathcal C}:\H\to R$ such that
$$
f_{\mathcal C}(T_w)=
\begin{cases}
1&\text{ if $w\in \mathcal C$}\\
0&\text{ otherwise}
\end{cases}
$$
whenever $w\in W$ has minimal length in its conjugacy class.
\end{thm}
Suppose $\mathcal C$ is a conjugacy class of the symmetric group.
An element $x\in \mathcal C$ has minimal length iff it is a product of distinct simple reflections.
If $x$ and $y$ are two such elements, the corresponding braid group elements are also conjugate (Tits);
hence $T_x$ and $T_y$ are conjugate in the Hecke algebra and $\chi(T_x)=\chi(T_y)$ for any trace function $\chi$.
In \cite[Thm 1.3]{Ram}, Ram shows that a trace function is determined by its value on such elements.

Consider the transition matrix $M=M(n)$ which expresses symmetric polynomials in the Jucys-Murphy operators
as linear combinations of the Geck-Rouquier basis of the centre.
If $\mathcal C$ is a conjugacy class of $\mathfrak S_n$ and $\lambda$ is a partition of $|\lambda |\le n$
let
$$
M_{\mathcal C,\lambda}:=\langle T_{w_{\mathcal C}},m_\lambda(L_1,\dots,L_n)\rangle
$$
where $w_{\mathcal C}\in \mathcal C$ is a fixed representative of minimal length
and $m_\lambda$ is the monomial symmetric function \cite[Section I.2]{Macdonald}.
This matrix is block upper triangular in view of Corollary~\ref{thm:alt} (Mathas).
The diagonal blocks are the submatrices $M^{(k)}=M^{(k)}(n)$ $(k\le n)$
indexed by the conjugacy classes $\mathcal C$ such that $\ell(w_{\mathcal C})=k$ and
partitions $\lambda$ such that $|\lambda |=k$.
The Dipper-James conjecture is equivalent to the following result for all diagonal blocks:
\begin{thm}\label{thm:altdj}
The columns of $M^{(k)}(n)$ span all of $R^d$ where $d$ is the number of rows.
\end{thm}

We prove this theorem by first establishing a special case conjectured by James.
If $n\ge 2k$, then the matrix $M^{(k)}(n)$ is square and independent of $n$ \cite[Thm 3.2]{Mat99}.
\begin{con}[James]\label{thm:jamesconj}\cite[Conj 3.5]{Mat99}
$M^{(k)}(2k)$ is invertible over $R$.
\end{con}
Given James' conjecture, Mathas argues in \cite[Thm 3.6]{Mat99} that
the centre has a basis consisting of a specific family of symmetric polynomials in $L_1,\dots,L_n$.
However we do not follow the last stage of the proof.

Mathas conjectures \cite[Conj 3.7]{Mat99} an explicit inverse for the matrix $M^{(k)}(2k)$.
We were inspired by this conjecture to study an analogue of $m_\lambda(L_1,\dots,L_n)$ indexed by compositions
instead of partitions. This idea plays a key role in our proof of James' conjecture.
Unfortunately the analogue of Mathas' conjecture for compositions turned out to be false when $k=5$.
Mathas' conjecture remains open.
Nevertheless, we do exhibit a formula for the inverse.

Another classical link between the symmetric functions in Jucys-Murphy elements
and the centre of the group algebra of the symmetric group is a result of Farahat and Higman~\cite{FH59}.
They derive a formula for elementary symmetric functions of classical Jucys-Murphy elements in terms of class sums,
and show that these generate the whole centre.
An analogous formula for the elementary symmetric functions in Jucys-Murphy elements holds in the Iwahori-Hecke algebra,
and consequently a corresponding set of generators for the centre of the Iwahori-Hecke algebra
can be obtained (Corollary~\ref{cor:gens.for.Z}).

The cyclotomic Hecke algebra and the Hecke algebra of type $B$ are also quotients of $\widehat \H$.
The reader might be tempted to consider the question of whether these quotients also preserve
the centre (assuming some invertible elements such as $q-1$).
However this is known to be false: Ariki has found a counter-example in the Hecke algebra of $G(3,1,2)$~\cite[Section 2]{Ariki96}.

The outline of the paper is as follows:
\begin{enumerate}\setcounter{enumi}{1}
\item Fix notation for compositions and sketch the connection with finite, totally ordered sets.
\item Develop properties of the space of quasi-symmetric polynomials, an analogue of symmetric polynomials indexed by compositions instead of partitions.  They restrict from $n+1$ variables to $n$ variables in a simple manner.
\item Study certain polynomials $\bm{\mathfrak a}(n)$ \cite[Def 2.17]{Mat99} of Mathas
which arise as coefficients in the powers of Hecke algebra generators.
\item Define some matrices that we use to establish James' conjecture.
\item Establish the link with the Hecke algebra. This generalises Proposition 2.21 of \cite{Mat99},
giving the coefficient of any increasing element $T_w$ of $\H_n$ in a product of Jucys-Murphy elements.
\item Prove the James and the Dipper-James Conjectures.
\end{enumerate}

We thank Alain Lascoux for his comments on an earlier version of this paper.
We dedicate this paper to Professor Gordon James.
The second author in particular thanks him for his support and encouragement.

\section{Compositions}\label{sec:comps}

\begin{defn} A \emph{composition} of $n$ is a sequence $\lambda=(\lambda_1,\dots,\lambda_l)$
of positive integers such that $\sum_{i}\lambda_i=n$.
In this case write $\ell(\lambda)=l$ and $|\lambda|=n$.
There is a unique composition $\emptyset$ of $0$.
\end{defn}

Let $\Lambda_n$ denote the set of compositions of $n$
and $\Lambda=\bigcup_{n}\Lambda_n$ denote the set of all compositions.
If $n>0$, there is a bijection between compositions of $n$ and compositions of smaller size
which takes $\lambda=(\lambda_1,\dots,\lambda_l)$ to $\lambda':=(\lambda_1,\dots,\lambda_{l-1})$.

In order to write down matrices indexed by compositions, it is convenient to list them is a fixed order.
We define a (listing) order $\Lambda$ recursively as follows:
For any pair $\lambda$ and $\mu$ in $\Lambda$,
\begin{equation}\label{eq:def.order}
\lambda<\mu\text{ iff }
\begin{cases}
|\lambda|<|\mu|,&\text{ or }\\
0<|\lambda|=|\mu| \text{ and } \lambda'<\mu'.
\end{cases}
\end{equation}

If $\lambda$ is a composition, let $\hat\lambda$ denote the partition with the same parts.
Among the set of compositions with the same set of parts, the partition is the last one listed.

In this paper, we make extensive use of total preorders.
The compositions of $n$ arise as the quotients of a totally ordered set of cardinality $n$
in the category of monotone functions.
A \emph{preorder} on a set $P$ is a relation $\preceq$ such that
    \begin{enumerate}
      \item\label{def:preorder:transitive} $x\preceq y$ and $y\preceq z$ implies $x\preceq z$ ($x,y,z\in P$), and
      \item\label{def:preorder:reflexive} $x\preceq x$ ($x\in P$).
    \end{enumerate}
The preorder $\preceq$ on $P$ is an \emph{order} if
    \begin{enumerate}\setcounter{enumi}{2}
    \item $x\preceq y$ and $y\preceq x$ imply $x=y$ ($x,y\in P$).
    \end{enumerate}
The preorder $\preceq$ on $P$ is \emph{total} if
    \begin{enumerate}\setcounter{enumi}{3}
      \item $x\preceq y$ or $y\preceq x$ ($x,y\in P$).
    \end{enumerate}
A \emph{totally ordered set} is a pair $\ul P=(P,\preceq_P)$ where
$P$ is a set and $\preceq_P$ is a total order on $P$.

If $\preceq_P$ and $\preceq_Q$ are preorders on sets $P$ and $Q$, a function $f:P\to Q$ is called \emph{monotone} if
$x\preceq_P y\implies f(x)\preceq_Q f(y)$ ($x,y\in P$).
(We reserve the term \emph{order-preserving} for when the converse is also true.)
A \emph{homomorphism of totally ordered sets} is a monotone function between the underlying sets.

The set $\ul n=\{1,2,\dots,n\}$ with the usual order $\le$ is a totally ordered set.
Every totally ordered set of cardinality $n$ is isomorphic to $\ul n$,
and the isomorphism is unique.

Let $P$ be a set with a preorder $\preceq$.
Define a relation $\sim$ on $P$ by
$x\sim y \iff x\preceq y \text{ and } y\preceq x$ $(x,y\in P)$.
Then $\sim$ is an equivalence relation.
The number $\ell(\preceq)$ of equivalence classes in $P$ is called the \emph{length} of $\preceq$.
If $x\in P$, let $[x]$ denote the equivalence class containing $x$.
The set $[P]$ of equivalence classes inherits an order given by
$[x]\preceq [y] \iff x\preceq y$ for $x,y\in P$.
Let $q_{\preceq}$ denote the function $P\to [P]$ which takes $x$ to $[x]$.

Let $\ul P=(P,\preceq_P)$ be a finite totally ordered set.
A \emph{composition} of $\ul P$ is a preorder $\preceq$ on $P$
such that $x\preceq_P y\implies x\preceq y$ for $x,y\in P$.
This is equivalent to the function $q_\preceq:\ul P\to ([P],\preceq)$ being monotone.
Conversely, a monotone function $f:\ul P\to \ul Q$ induces a composition $\preceq_f$ of $\ul P$
given by $x\preceq_f y\iff f(x)\preceq_Q f(y)$ ($x,y\in P$).
Given compositions $\preceq$ and $\preceq'$ of $\ul P$, we say
$\preceq$ \emph{is contained by} $\preceq'$ (denoted $\preceq\subseteq\preceq'$)
if $x\preceq y\Leftarrow x\preceq' y$ $(x,y\in P)$.
This is the case iff $q_{\preceq}$ factors through $q_{\preceq'}$.

Suppose $P$ above has cardinality $n$.
There is a one to one correspondence between compositions of $P$ and $n$ as follows.
Let $\preceq$ be a composition of $P$ of length $k=\ell(\preceq)$ and let
$\theta:([P],\preceq)\to\ul k$ be the unique isomorphism.
For $i\le k$, let $\lambda_i$ be the cardinality of the inverse image of $i$ under $\theta\circ q_\preceq$.
Then $\lambda_\preceq=(\lambda_1,\dots,\lambda_k)$ is a composition of $n$.
We call $\lambda$ the \emph{shape} of $\preceq$.
Conversely, if $\lambda$ is a composition, we denote the corresponding composition of $\ul n$ by $\preceq_\lambda$.
For example, the composition $(2,3)$ of $5$ corresponds the total preorder $\preceq_{(2,3)}$ of $\{1,2,3,4,5\}$
where $1$ and $2$ are smaller than all elements, and $3$, $4$ and $5$ are larger than all elements.

\section{Quasi-symmetric Polynomials}\label{sec:polynomials}

Let $R$ be a commutative ring with 1 and $R[X_1,\dots,X_n]$ denote the ring of polynomials over $R$ in $n$ independent variables.
In this section we study the $R$-subalgebra of \emph{quasi-symmetric polynomials}.
These polynomials are similar to symmetric polynomials, but have a basis indexed by compositions with at most $n$ parts, rather than by partitions.

\begin{defn}
If $\preceq$ is a total preorder on a set $P$, the polynomial
$$
p^\preceq_{\le n}=p^\preceq(X_1,\dots,X_n)=\sum_{\substack{f:(P,\preceq)\to \ul n\\ \text{order-preserving}}}X^f
\quad\text{where}\quad
X^f=\prod_{i\in P} X_{f(i)}.
$$
is called a monomial quasi-symmetric polynomial.
The span $QSym_{\le n}$
of these polynomials is called the set of \emph{quasi-symmetric polynomials}.
\end{defn}

For example,
$$
p^{(1,2)}(X_1,\dots,X_4)=X_1X_2^2+X_1X_3^2+X_1X_4^2+X_2X_3^2+X_2X_4^2+X_3X_4^2\\
$$

The study of quasi-symmetric polynomials dates at least as far back as \cite{Gessel84}, and has more recently had further attention through for instance~\cite{Stanley95} and~\cite{GelfandEtal95}.

The following Lemma is well known.
\begin{lem}
The family $p^\lambda_{\le n}=p^{\preceq_\lambda}_{\le n}$ indexed by compositions $\lambda$ with at most $n$ parts is a basis of the space $QSym_{\le n}$ of quasi-symmetric polynomials.
\end{lem}
\begin{proof}
Let $\preceq$ be a total preorder on a set $P$ of cardinality $k$.
We say $x$ is a \emph{minimum} if $x\preceq y$ for all $y\in P$. (There may be more than one).
By repeatedly selecting a minimum, we list the elements $p_1,p_2,\dots,p_k$ of $P$ is increasing order.
The resulting bijective monotone function $p:\ul k\to P$ induces a composition $\preceq_k$ of $\ul k$.
We have $p^{\preceq_k}_{\le n}=p^\preceq$.
Hence every monomial quasi-symmetric polynomial has the form $p^\lambda_{\le n}$ for some composition $\lambda$.

If $\preceq$ is a total preorder on a set $P$,
there exists an order preserving $f:P\to \ul n$ iff $\ell(\preceq_P)\le n$.
Hence $p^\lambda_{\le n}$ is non-zero iff $\lambda$ has at most $n$ parts.

It remains to show that the family are linearly independent.
Consider a monomial $X_1^{i_1}\cdots X_n^{i_n}$ of degree $k$.
The monomial has the form $X^f$ for some monotone function $f:\ul k\to \ul n$,
which induces a composition $\preceq_f$ of $\ul k$.
Although $f$ is not unique, the resulting composition is; the monomial contributes only to $p^{\preceq_f}$.
\end{proof}

$QSym_{\le n}$ is a subalgebra of $R[X_1,\dots,X_n]$ thanks to the following:

\begin{prop}\label{thm:pmult}
If $P$ and $Q$ are disjoint, finite sets with total preorders $\preceq_P$ and $\preceq_Q$ respectively, then as polynomials of $X_1,\dots,X_n$,
\[
p^{\preceq_P} p^{\preceq_Q} =\sum_{\preceq}p^{\preceq}
\]
where the sum varies over preorders $\preceq$ of $P\cup Q$ which restrict to $\preceq_P$ on $P$ and $\preceq_Q$ on $Q$.
\end{prop}
\begin{proof}
The terms on the left hand side are indexed by pairs of order preserving functions $f:P\to \ul n$ and $g:Q\to \ul n$.
The union $h=f\cup g:P\cup Q\to \ul n$ induces a preorder $\preceq$ on $P\cup Q$ by
$x\preceq y\iff h(x)\le h(y)$. These index terms on the right hand side.
\end{proof}

The space of quasi-symmetric polynomials has another interesting basis.
\begin{defn}
If $\preceq$ is a total preorder on a set $P$, define
$$
q^\preceq_{\le n}=q^\preceq(X_1,\dots,X_n)=(-1)^{\ell(\preceq)}\sum_{\substack{f:(P,\preceq)\to \ul n\\ \text{monotone}}}X^f.
$$
\end{defn}
For example,
$$
q^{(3,2)}(X_1,X_2,X_3)=X_1^3X_2^2+X_1^3X_3^2+X_2^3X_3^2+X_1^5+X_2^5+X_3^5.
$$

\begin{lem}\label{lem:inversion}
If $\preceq$ is a composition of finite, totally ordered set $\ul P$, then
\begin{align*}
(-1)^{\ell(\preceq)} q^{\preceq}&=\sum_{\preceq'\subseteq\preceq} p^{\preceq'}\quad\text{and}\\
(-1)^{\ell(\preceq)} p^{\preceq}&=\sum_{\preceq'\subseteq\preceq} q^{\preceq'}
\end{align*}
as polynomials of $X_1,\dots,X_n$,
where the sums vary over compositions $\preceq'$ of $\ul P$ contained in $\preceq$.
\end{lem}

\begin{proof}
The first equation says
$$
\sum_{\substack{f:(P,\preceq)\to \ul n\\ \text{monotone}}}X^f
=\sum_{\preceq'\subseteq\preceq}\sum_{\substack{f:(P,\preceq')\to \ul n\\ \text{order-preserving}}}X^f.
$$
The terms on the left hand side are indexed by monotone $f:(P,\preceq)\to \ul n$.
Each such function induces a composition $\preceq'$ of $\ul P$ by
$x\preceq' y\iff f(x)\le f(y)$.
Note that $\preceq$ contains $\preceq'$ and that $f:(P,\preceq')\to \ul n$ is order preserving.
Such pairs index the right hand side.

It remains to verify the second equation.
Let $P'$ denote $P$ with its maximum removed.
There is a bijection $I$ between the set of compositions of $\ul P$
and the power set of $P'$ which preserves the meaning of ``contains''.
If $\preceq$ is a composition, let $I(\preceq)$ be the set of $x\in P'$ such
$y\preceq x\implies y\preceq_P x$ $(y\in P)$
The set $I(\preceq)$ has cardinality $\ell(\preceq)-1$.
The M\"obius function of the power set is well known \cite[3.8.3]{Stanley} to be
$\mu(\preceq, \preceq')=(-1)^{\ell(\preceq)-\ell(\preceq')}$.
The second equation is the M\"obius inversion formula \cite[3.7.1]{Stanley} (or the inclusion-exclusion principle).
\end{proof}

\begin{prop} 
If $P$ and $Q$ are disjoint, finite sets with total preorders $\preceq_P$ and $\preceq_Q$ respectively,
then
\[
\sum_{\preceq}(-1)^{\ell(\preceq)}=(-1)^{\ell(\preceq_P)+\ell(\preceq_Q)}.
\]
where the sum varies over preorders $\preceq$ of $P\cup Q$ which restrict to $\preceq_P$ on $P$ and $\preceq_Q$ on $Q$.
\end{prop}

\begin{rem}This is equivalent to the multinomial identity
\[\sum_{i\le\min(a,b)}(-1)^i\binom{a+b-i}{i,a-i,b-i}=1\]
for non-negative integers $a$ and $b$.
\end{rem}

\begin{proof}
We prove this result by induction on $|P|+|Q|$.

If $P=\emptyset$ or $Q=\emptyset$ then the statement is clear.
Suppose then that $P$ and $Q$ are non-empty and the statement holds if $P\cup Q$ is smaller.
Let $\mathcal M$ (resp. $\mathcal N$) be the set of maximum elements of $P$ under $\preceq_P$ (resp. $Q$ under $\preceq_Q$).

Let $\preceq$ be a preorder on $P\cup Q$ which restricts to $\preceq_P$ on $P$ and $\preceq_Q$ on $Q$,
and consider the set $\mathcal U$ of maximal elements in $(P\cup Q,\preceq)$.
Precisely one of the following is true:
    \begin{enumerate}\renewcommand{\labelenumi}{(\roman{enumi})}
    \item $\mathcal U=\mathcal M$,
    \item $\mathcal U=\mathcal N$, or
    \item $\mathcal U=\mathcal M\cup\mathcal N$.
    \end{enumerate}
Let $\preceq_P'$ be the restriction of $\preceq_P$ to $P'=(P-\mathcal M)$,
define $\preceq_Q'$ and $\preceq'$ similarly.
It follows that the sum on the left hand side decomposes into three parts:
\begin{align*}
\sum_{\preceq}  (-1)^{\ell(\preceq)}
 &=(-1)\sum_{\preceq'\text{ on } P'\cup  Q}  (-1)^{\ell(\preceq')}\\
 &+(-1)\sum_{\preceq'\text{ on } P\cup Q'}  (-1)^{\ell(\preceq')}\\
 &+(-1)\sum_{\preceq'\text{ on } P'\cup Q'}  (-1)^{\ell(\preceq')}
\end{align*}
Applying the inductive hypothesis, and the fact $\ell(\preceq_P)=\ell(\preceq_P')+1$,
\begin{align*}
 &=(-1)\Big( (-1)^{\ell(\preceq_P')+\ell(\preceq_Q)}+(-1)^{\ell(\preceq_P)+\ell(\preceq_Q')}\\
 &\qquad\qquad +(-1)^{\ell(\preceq_P')+\ell(\preceq_Q')}\Big)\\
 &=(-1)\left((-1)^{-1}+(-1)^{-1}+(-1)^{-2}\right)(-1)^{\ell(\preceq_P)+\ell(\preceq_Q)}\\
 &=(-1)^{\ell(\preceq_P)+\ell(\preceq_Q)}.
\end{align*}
\end{proof}

\begin{prop}\label{thm:sumsP&QwithM}
If $P$ and $Q$ are disjoint, finite sets with total preorders $\preceq_P$ and $\preceq_Q$ respectively,
then as polynomials of $X_1,\dots,X_n$,
\[
q^{\preceq_P} q^{\preceq_Q} =\sum_{\preceq}q^{\preceq}
\]
where the sum varies over preorders $\preceq$ of $P\cup Q$ which restrict to $\preceq_P$ on $P$ and $\preceq_Q$ on $Q$.
\end{prop}

\begin{proof}
The statement for $n=0$ is vacuous; the statement for $n=1$ is equivalent to the previous proposition.

Choose monotone $f:P\to\ul n$ and $g:Q\to\ul n$.
Let $P_i=f^{-1}(i)$, $Q_i=g^{-1}(i)$ and note that $P_i\cup Q_i$ is the inverse image of $i$
under $h=f\cup g:P\cup G\to\ul n$.
Let ${\preceq_P}_i$ denote the restriction of $\preceq_P$ to $P_i$, and similarly for ${\preceq_Q}_i$.

There is a bijection between preorders $\preceq$ on $P\cup Q$ such that $h:(P\cup Q,\preceq)\to\ul n$ is monotone,
and $n$-tuples of preorders $\preceq_i$ on $P_i\cup Q_i$ given by restriction.
The restriction of $\preceq$ to $P$ agrees with $\preceq_P$ iff each corresponding $\preceq_i$ restricts to ${\preceq_P}_i$.
The following sums are indexed by such preorders:
\begin{align*}\sum_{\substack{\preceq\\
                h:P\cup Q\to\ul n}}(-1)^{l(P\cup Q,\preceq)}
 &=\prod_{i\in\ul n}\left(\sum_{\preceq_i}(-1)^{l(\preceq_i)}\right)\\
 &=\prod_{i\in\ul n}(-1)^{l({\preceq_P}_i)+l({\preceq_Q}_i)}\\
 &=(-1)^{l(\preceq_P)+l(\preceq_Q)},
\end{align*}
where the second equality is a consequence of the previous proposition
and the third follows since $(-1)^{l(\preceq_P)}=\prod_{i\in\ul n}(-1)^{l({\preceq_P}_i)}$.

Summing over all pairs of homomorphisms $(f,g)$ yields
\begin{align*}
 \sum_{\preceq}q_{\le n}^{\preceq}
 &=\sum_{\substack{f:P\to\ul n\\ g:Q\to\ul n}}\ \sum_{\substack{\preceq \\
                 f\cup g:P\cup Q\to\ul n}}(-1)^{l(\preceq)}X^{f\cup g}\\
 &=\sum_{\substack{f:P\to\ul n\\ g:Q\to\ul n}}(-1)^{l(\preceq_P)+l(\preceq_Q)}X^fX^g\\
 &=q_{\le n}^{\preceq_P}q_{\le n}^{\preceq_Q}.
\end{align*}
\end{proof}

\section{The invertibility of certain power series}\label{sec:xy=1}

Let $R=\Z[q,q^{-1}]$.
If $n\in \Z$, let $[n]_q\in R$ denote the unique Laurent polynomial such that $q^n-q^{-n}=[n]_q (q-q^{-1})$.
Let $\xi=q^{-1}(q-1)^2$. Define
$$
\bm{\mathfrak a}(0):=1,\quad
\bm{\mathfrak a}(s):=[s]_q\xi,\quad
\bm{\mathfrak b}(0):=1,\quad\text{and}\quad
\bm{\mathfrak b}(s):=-s\xi,
$$
where $s>0$ is an integer.

\begin{prop}\label{prop:abinv}
$$
a(X)=\sum_{s\ge 0} \bm{\mathfrak a}(s)X^s\quad\text{and}\quad
b(X)=\sum_{s\ge 0} \bm{\mathfrak b}(s)X^s
$$
are inverse in the ring $R[[X]]$ of power series over $R$.
\end{prop}
\begin{proof}
\begin{align*}
\sum_{s\ge 0} \bm{\mathfrak a}(s)X^s
        &=1+\frac{\xi}{q-q^{-1}}\sum_{s\ge 0}(qX)^s-(q^{-1}X)^s\\
        &=1+\frac{\xi}{q-q^{-1}}\left(\frac{1}{1-qX}-\frac{1}{1-q^{-1}X}\right)\\
        &=\frac{(1-X)^2}{(1-qX)(1-q^{-1}X)}.\\
\sum_{s\ge 0} \bm{\mathfrak b}(s)X^s
        &=1-\xi X\sum_{s\ge 1}sX^{s-1}\\
        &=1-\xi X\frac{1}{(1-X)^2}\\
        &=\frac{(1-qX)(1-q^{-1}X)}{(1-X)^2}.
\end{align*}
\end{proof}

The following corollary shows that $\bm{\mathfrak a}$ (but not $\bm{\mathfrak b}$)
is the same as the one used by Mathas in \cite[Lemma 2.16(iii)]{Mat99}.
Since $\bm{\mathfrak b}(0)$ and $\bm{\mathfrak a}(0)$ are $1$,
this recurrence characterises both $\bm{\mathfrak a}$ and $\bm{\mathfrak b}$
in terms of the other.
\begin{cor}\label{lem:ab}
$$
\sum_{s+t=r}\bm{\mathfrak a}(s)\bm{\mathfrak b}(t)=
\begin{cases}
1&\text{if $r=0$,}\\
0&\text{otherwise.}
\end{cases}
$$
\end{cor}

Totally ordered sets of cardinality $1$ are the terminal objects---there exists a unique
monotone function $f:\ul P\to \ul 1$ from any totally ordered set $P$.
The induced preorder $\preceq_f$ of $\ul P$ is given by $x\preceq_f y$ for all $x,y\in P$.
This is called the \emph{trivial} composition.

\begin{prop}\label{prop:ab2comp}
If $\ul P$ is a totally ordered set of cardinality $n$, then
$$
\bm{\mathfrak a}(n)=\sum_{\substack{r:\ul P\to \ul P\\ r^2=r}}\xi^{\ell(\preceq_r)}
\quad\text{and}\quad
\bm{\mathfrak b}(n)=\sum_{\substack{r:\ul P\to \ul P\\ r^2=r\\ \preceq_r\text{ is trivial}}}(-\xi)^{\ell(\preceq_r)}.
$$
\end{prop}
\begin{proof}
The second formula is immediate. If $n=0$ and $P$ is empty,
there is a unique function $r:P\to P$ and it contributes $1$ to the sum.
Alternately, if $n>0$ and $P$ is non-empty, the whole of $P$ is mapped to one element so the sum yields $-n$.

We prove the first formula by induction on $n$.
It is trivial when $n=0$.
Suppose that $\ul P=\ul n$ is non-empty.
Any idempotent monotone function partitions $\ul n$ into the set $\mathcal M$ of maximum elements relative to $\preceq_r$
and a complement of the form $\ul t$ for some $t<n$.
If we restrict the idempotent to $\mathcal M$ of cardinality $s=n-t$, we obtain an idempotent inducing the trivial order.
If we restrict the idempotent to $\ul t$ we obtain another (arbitrary) idempotent.
Hence
\begin{align*}
\sum_{\substack{r:\ul n\to \ul n\\ r^2=r}}\xi^{\ell(\preceq_r)}
&=\sum_{\substack{s+t=n\\ s>0}} \Big(\sum_{\substack{r:\mathcal M\to \mathcal M\\ r^2=r\\ \preceq_r\text{ is trivial}}}\xi^{\ell(\preceq_r)}\Big)
\Big(\sum_{\substack{r:\ul t\to \ul t\\ r^2=r}}\xi^{\ell(\preceq_r)-1}\Big)\\
&=-\sum_{\substack{s+t=n\\ s>0}} \bm{\mathfrak b}(s) \bm{\mathfrak a}(t).
\end{align*}
The result now follows from Corollary~\ref{lem:ab}.
\end{proof}

Mathas \cite[Def 2.17]{Mat99} uses the following formula to define $\mathfrak a$.
\begin{cor}\label{cor:reltomathas}
If $n$ is a positive integer,
$$
\bm{\mathfrak a}(n)=\sum_{m=1}^n \binom{n+m-1}{2m-1} \xi^m.
$$
\end{cor}
\begin{proof}
The idempotent monotone functions $r:\ul n\to \ul n$ with $m$ fixed points may be enumerated as follows:
Take $n+m+1$ boxes and arrange them in a row. Label the first $j_0$ and the last $j_m$.
Choose any $2m-1$ from the remaining $n+m-1$ boxes, and label them alternately $i_1$, $j_1$, $i_2$, $j_2$,\dots,$j_{m-1}$,$i_m$.
Now place the numbers $1$, $2$,\dots,$n$ into the boxes which are not labelled with a $j$.
The corresponding function $r:\ul n\to \ul n$
maps the numbers in boxes between $j_{k-1}$ and $j_k$ to the number in box $i_k$.
\end{proof}

If $\lambda$ is a composition with $l$ parts, define
$$
\bm{\mathfrak a}(\lambda)=\prod_{1\le i\le l}\bm{\mathfrak a}(\lambda_i)
\quad\text{and}\quad
\bm{\mathfrak b}(\lambda)=\prod_{1\le i\le l}\bm{\mathfrak b}(\lambda_i).
$$

Examples of $\bm{\mathfrak a}(\lambda)$ and $\bm{\mathfrak b}(\lambda)$ for small $n$ are given in Table~\ref{tab:a-b-functions}.

\begin{table}[ht]
\caption{$\bm{\mathfrak a}(\lambda)$ and $\bm{\mathfrak b}(\lambda)$ for $\lambda\in\Lambda_n$, $n\le 3$.}\label{tab:a-b-functions}
\begin{tabular}{llll}
$n$\quad\ &$\lambda$&$\bm{\mathfrak a}(\lambda)$&$\bm{\mathfrak b}(\lambda)$\\
\hline
0&$\emptyset$&$1$&$1$\\
1&$(1)$&$\xi$&$-\xi$\\
2&$(2)$&$\xi^2+2\xi$&$-2\xi$\\
&$(1,1)$&$\xi^2$&$\xi^2$\\
3&$(3)$&$\xi^3+4\xi^2+3\xi$&$-3\xi$\\
&$(1,2)$&$\xi^3+2\xi^2$&$2\xi^2$\\
&$(2,1)$&$\xi^3+2\xi^2$&$2\xi^2$\\
&$(1,1,1)$&$\xi^3$&$-\xi^3$\\
\end{tabular}
\end{table}

\begin{lem}\label{cor:abrel}
If $\lambda$ is a composition,
$$
\bm{\mathfrak a}(\lambda)=\sum_{\mu\supseteq\lambda}(-1)^{\ell(\mu)}\bm{\mathfrak b}(\mu)
\quad\text{and}\quad
\bm{\mathfrak b}(\lambda)=\sum_{\mu\supseteq\lambda}(-1)^{\ell(\mu)}\bm{\mathfrak a}(\mu).
$$
\end{lem}
\begin{proof}
Suppose $\ul P$ is a totally ordered set of cardinality $n$ and $\preceq$ be a composition of shape $\lambda$.
Then the equivalences classes $P_1,\dots,P_l$ of $P$ have size $\lambda_i=|P_i |$.
An idempotent such that $\preceq\subseteq\preceq_r$
corresponds to a family $r_i$ of idempotents on the classes $P_i$.
Applying Proposition~\ref{prop:ab2comp}, we find:
\begin{align*}
\bm{\mathfrak a}(\lambda)
=\prod_{1\le i\le l}\Big(\sum_{\substack{r_i:\ul P_i\to \ul P_i\\ r_i^2=r_i}}\xi^{\ell(\preceq_{r_i})}\Big)
=\sum_{\substack{r:\ul P\to \ul P\\ r^2=r\\ \preceq\subseteq\preceq_r}}\xi^{\ell(\preceq_r)}.\\
\bm{\mathfrak b}(\lambda)
=\prod_{1\le i\le l}\Big(-\sum_{\substack{r_i:\ul P_i\to \ul P_i\\ r_i^2=r_i\\ \preceq_{r_i}\text{ is trivial}}}\xi^{\ell(\preceq_{r_i})}\Big)
=(-1)^{\ell(\preceq)}\sum_{\substack{r:\ul P\to \ul P\\ r^2=r\\ \preceq=\preceq_r}}\xi^{\ell(\preceq_r)}.
\end{align*}
This proves the first statement. The second one follows by M\"obius Inversion.
\end{proof}

Consider the power series
$$
a(X_1,\dots,X_n)=\prod_{1\le i\le n} a(X_i)
$$
where $a(X)$ is as defined in Proposition~\ref{prop:abinv},
and denote the homogeneous component of degree $k$ by $a_k(X_1,\dots,X_n)$.

\begin{cor}\label{cor:basisdecomp}
\begin{align*}
a_k(X_1,\dots,X_n)&=\sum_{|\lambda |=k} \bm{\mathfrak a}(\lambda)p^\lambda_{\le n}=\sum_{|\lambda |=k} \bm{\mathfrak b}(\lambda)q^\lambda_{\le n},\\
b_k(X_1,\dots,X_n)&=\sum_{|\lambda |=k} \bm{\mathfrak b}(\lambda)p^\lambda_{\le n}=\sum_{|\lambda |=k} \bm{\mathfrak a}(\lambda)q^\lambda_{\le n}.
\end{align*}
\end{cor}

We require the following technical result for the proof of Proposition~\ref{thm:reduction1} below.
\begin{lem}\label{lem:recurrence} If $n>0$,
$$
\sum_{0\le k<r} a_k(X_1,\dots,X_{n-1})=\sum_{0\le s<r} \bm{\mathfrak b}(s)X_n^s\sum_{0\le t<r-s} a_t(X_1,\dots,X_n).
$$
\end{lem}
\begin{proof}
Since $a(X_n)$ and $b(X_n)$ are inverse, we have
$a(X_1,\dots,X_{n-1})=a(X_1,\dots,X_{n-1})a(X_n)b(X_n)=a(X_1,\dots,X_n)b(X_n)$.
Comparing terms of degree $k=s+t$ less than $r$ yields the recurrence.
\end{proof}

\section{Matrices}

Our goal is to calculate $\l T_w, p^\mu(\LL_1,\dots,\LL_n) \r$ for
increasing $w\in W$ and compositions $\mu$ such that $\ell(w)=|\mu |$.
This bilinear form is independent of several choices,
but this is only apparent to us because they satisfy the same recurrence.
We introduce this recurrence by means of certain square matrices
indexed by compositions of size less than $k$.

If $\lambda$ is a composition with $l$ parts and $k\le l$,
we call the composition $\mu=(\lambda_1,\lambda_2,\dots,\lambda_k)$ a \emph{prefix} of $\lambda$.
Recall from Section~\ref{sec:comps} that for $\lambda=(\lambda_1,\lambda_2,\dots,\lambda_l)$ a composition of $n$, we set $\lambda'=(\lambda_1,\lambda_2,\dots,\lambda_{l-1})$; $\lambda'$ is a particular prefix of $\lambda$.

\begin{defn}
If $\lambda$ and $\mu$ are compositions, define
\begin{align*}
J_{\lambda,\mu}&=
\begin{cases}
(-1)^{\ell(\mu)}&\text{if $\lambda\subseteq \mu$,}\\
0\qquad\qquad\quad&\text{otherwise,}
\end{cases}\\
K_{\lambda,\mu}&=
\begin{cases}
(-1)^{\ell(\mu)}&\text{if $\lambda\subseteq \nu$ for some prefix $\nu$ of $\mu$,}\\
0\qquad\qquad\quad&\text{otherwise,}
\end{cases}\\
Z_{\lambda,\mu}&=
\begin{cases}
1&\text{if $\lambda= \mu$ or $\lambda= \mu'$,}\\
0\qquad\qquad\quad&\text{otherwise, and}
\end{cases}\\
Y_{\lambda,\mu}&=
\begin{cases}
(-1)^{\ell(\mu)-\ell(\lambda)}&\text{if $\lambda$ is a prefix of $\mu$,}\\
0\qquad\qquad\quad&\text{otherwise.}
\end{cases}
\end{align*}
\end{defn}

Let $k$ be a positive integer,
and consider the matrices $J=J^{(k)}$, $K=K^{(k)}$, $Z=Z^{(k)}$ and $Y=Y^{(k)}$
indexed by compositions $\lambda$ such that $|\lambda |<k$, listed in the order specified in Section~\ref{sec:comps}.

For example, for $k=4$ (compositions of 3 or less) we have:

\[\tiny
J=
\left(
\begin{array}{rrrrrrrr}
1&&&&&&&\\
&-1&&&&&&\\
&&-1&1&&&&\\
&&&1&&&&\\
&&&&-1&1&1&-1\\
&&&&&1&&-1\\
&&&&&&1&-1\\
&&&&&&&-1
\end{array}\right)\
K=\left(
\begin{array}{rrrrrrrr}
1&-1&-1&1&-1&1&1&-1\\
&-1&&1&&1&&-1\\
&&-1&1&&&1&-1\\
&&&1&&&&-1\\
&&&&-1&1&1&-1\\
&&&&&1&&-1\\
&&&&&&1&-1\\
&&&&&&&-1
\end{array}\right)
\]
and
\[\tiny
Z=\left(
\begin{array}{rrrrrrrr}
1&1&1&&1&&&\\
&1&&1&&1&&\\
&&1&&&&1&\\
&&&1&&&&1\\
&&&&1&&&\\
&&&&&1&&\\
&&&&&&1&\\
&&&&&&&1
\end{array}\right)\
Y=\left(
\begin{array}{rrrrrrrr}
1&-1&-1&1&-1&1&1&-1\\
&1&&-1&&-1&&1\\
&&1&&&&-1&\\
&&&1&&&&-1\\
&&&&1&&&\\
&&&&&1&&\\
&&&&&&1&\\
&&&&&&&1
\end{array}\right).
\]

\begin{lem}\label{matrix1}
These matrices satisfy
$$
J^2=I,\quad
K^2=I,\quad
ZY=I,\quad
JK=Y\quad\text{and}\quad
KJ=Z.
$$
\end{lem}

\begin{proof}
Recall that there is a bijection between compositions of $k$ and smaller compositions.
The matrices for $k+1$ may be given recursively in terms of the matrices for $k$ as follows:
\begin{align*}
J^{(k+1)}&=
\begin{pmatrix}J&0\\ 0&-K\end{pmatrix},
\\
K^{(k+1)}&=
\begin{pmatrix}K&-K\\ 0&-K\end{pmatrix},
\\
Z^{(k+1)}&=
\begin{pmatrix}Z&I\\ 0&I\end{pmatrix},
\\
Y^{(k+1)}&=
\begin{pmatrix}Y&-Y\\ 0&I\end{pmatrix}.
\end{align*}
The equations can now be established by induction on $k$.
\end{proof}

\begin{defn}
If $\lambda$ and $\mu$ are compositions of $k$ and $l$, choose $n\ge\ell(\lambda)$
and define elements $A_{\lambda,\mu}$ and $B_{\lambda,\mu}$ of $R$ as follows:
If $k\ge l$,
\begin{align*}
a_{k-l}(X_1\cdots X_n) p^\mu_{\le n}&=\sum_{\ell(\nu)\le n} A_{\nu,\mu} p^\nu_{\le n}\quad\text{and}\\
b_{k-l}(X_1\cdots X_n) p^\mu_{\le n}&=\sum_{\ell(\nu)\le n} B_{\nu,\mu} p^\nu_{\le n}.
\end{align*}
If $k<l$, then $A_{\lambda,\mu}=0$ and $B_{\lambda,\mu}=0$.
\end{defn}

First we show that the coefficients are well defined.
\begin{lem}\label{lem:omnibus}
\

\begin{enumerate}
\item
If $\mu$ is a composition of $l$,
$a_{k-l}(X_1\cdots X_n) p^\mu(X_1,\dots,X_n)$ and $b_{k-l}(X_1\cdots X_n) p^\mu(X_1,\dots,X_n)$
are quasi-symmetric polynomials of degree $k$.
\item
$A_{\lambda,\mu}$ and $B_{\lambda,\mu}$ do not depend on the choice of $n$.
\end{enumerate}
\end{lem}
\begin{proof}
First we show (1). Fix $n$. Since $a_{k-l}(X_1\cdots X_n)$ is a symmetric polynomial in $X_1,\dots,X_n$, it is also a quasi-symmetric polynomial.
Since the product of quasi-symmetric polynomials is quasi-symmetric,
it follows that $a_{k-l}(X_1\cdots X_n) p^\mu_{\le n}$ is quasi-symmetric.
The family $\{p^\lambda_{\le n}\}$ is a basis for quasi-symmetric polynomials of degree $k$, so there exist unique coefficients $A_{\lambda,\mu}$ (with $n$ fixed).

Now consider (2).
If $\ul P$ and $\ul Q$ are disjoint totally ordered sets with compositions $\preceq_P$ and $\preceq_Q$
of type $\eta$ and $\mu$ respectively,
let $c^\lambda_{\eta,\mu}$ denote the number of $\preceq$ of type $\lambda$ indexing the sum in Proposition~\ref{thm:pmult}.
Using Corollary~\ref{cor:basisdecomp}, we have
$$
A_{\lambda,\mu}=\sum_{|\eta |=k} \bm{\mathfrak a}(\eta)c^\lambda_{\eta,\mu}.
$$
This does not depend on $n$.

The proof of (1) and (2) for $B_{\lambda,\mu}$ is analogous.
\end{proof}

Let $k$ be a positive integer,
and consider the matrices $A=A^{(k)}$ and $B=B^{(k)}$
indexed by compositions $\lambda$ such that $|\lambda |<k$, listed in the order specified in Section~\ref{sec:comps}.  For example, for $k=4$ we have

\def\aa{\bm{\mathfrak a}}
\def\bb{\bm{\mathfrak b}}

{\footnotesize
\begin{align*}
A&=
\left(
\begin{array}{cccccccc}
\aa(0)&&&&&&&\\
\aa(1)&\aa(0)&&&&&&\\
\aa(2)&\aa(1)&\aa(0)&&&&&\\
\aa(1)^2&2\aa(1)&&\aa(0)&&&&\\
\aa(3)&\aa(2)&\aa(1)&&\aa(0)&&&\\
\aa(1)\aa(2)&\aa(2)+\aa(1)^2&\aa(1)&\aa(1)&&\aa(0)&&\\
\aa(2)\aa(1)&\aa(2)+\aa(1)^2&\aa(1)&\aa(1)&&&\aa(0)&\\
\aa(1)^3&3\aa(1)^2&&3\aa(1)&&&&\aa(0)
\end{array}\right)\\
&=
\left(
\begin{array}{cccccccc}
1&&&&&&&\\
\xi&1&&&&&&\\
\xi^2+2\xi&\xi&1&&&&&\\
\xi^2&2\xi&&1&&&&\\
\xi^3+4\xi^2+3\xi&\xi^2+2\xi&\xi&&1&&&\\
\xi^3+2\xi^2&2\xi^2+2\xi&\xi&\xi&&1&&\\
\xi^3+2\xi^2&2\xi^2+2\xi&\xi&\xi&&&1&\\
\xi^3&3\xi^2&&3\xi&&&&1
\end{array}\right)\\
\intertext{\normalsize and similarly}
B&=
\left(
\begin{array}{cccccccc}
1&&&&&&&\\
-\xi&1&&&&&&\\
-2\xi&-\xi&1&&&&&\\
\xi^2&-2\xi&&1&&&&\\
-3\xi&-2\xi&-\xi&&1&&&\\
2\xi^2&-2\xi+\xi^2&-\xi&-\xi&&1&&\\
2\xi^2&-2\xi+\xi^2&-\xi&-\xi&&&1&\\
-\xi^3&3\xi^2&&-3\xi&&&&1
\end{array}\right).
\end{align*}}

\begin{lem}\label{lem:prodab}
$AB=I$.
\end{lem}
\begin{proof}
Let $k>0$ and choose $n=k$.
Let $\lambda$ and $\mu$ be compositions such that $|\lambda |,|\mu |<k$.
If $\alpha$ is a composition with $|\alpha |<k$,
we have $A_{\lambda,\alpha}B_{\alpha,\mu}=0$ unless $|\lambda |\ge |\alpha |\ge |\mu|$.
Hence
\begin{equation}
(AB)_{\lambda,\mu}=\sum_{|\alpha |<k} A_{\lambda,\alpha}B_{\alpha,\mu}=\sum_{|\lambda |\ge |\alpha |\ge |\mu|}A_{\lambda,\alpha}B_{\alpha,\mu}.\label{eq:sofar}
\end{equation}
Let $m=|\lambda |- |\mu|$.
If $m<0$, then \eqref{eq:sofar} vanishes, as required.
Assume then that $m\ge 0$, and consider the homogenous component of degree $m$ in the equation
$a(X_1,\dots,X_n)b(X_1,\dots,X_n)=1$ (Proposition \ref{prop:abinv}):
$$
\sum_{i+j=m} a_i(X_1,\dots,X_n)b_j(X_1,\dots,X_n)=
\begin{cases}
1&\text{if $m=0$,}\\
0&\text{otherwise}.
\end{cases}
$$
Multiply by $p^\mu$ and apply the definition of $A$ and $B$:
$$
\sum_{i+j=m} \sum_{|\alpha |=|\mu| + j\ }\sum_{|\beta |=|\mu |+m\ } A_{\beta,\alpha}B_{\alpha,\mu}p^\beta=
\begin{cases}
p^\mu&\text{if $m=0$,}\\
0&\text{otherwise}.
\end{cases}
$$
Taking the coefficient of $p^\lambda$, we find
$$
\sum_{|\lambda |\ge |\alpha |\ge |\mu|}A_{\lambda,\alpha}B_{\alpha,\mu}=
\begin{cases}
1&\text{if $\lambda=\mu$,}\\
0&\text{otherwise}.
\end{cases}
$$
\end{proof}

\begin{lem}\label{matrix2}
$AJ=JB$.
\end{lem}
\begin{proof}
Let $k>0$ and choose $n\ge k$.
The span $\mathfrak I_{\le n}$ of $\{ p^\eta_{\le n}\}$ indexed by compositions $\eta$ such that $|\eta |> n$
is an ideal of $QSym_{\le n}$.
The quotient $R$-algebra $QSym_{\le n}/\mathfrak I_{\le n}$ has two bases $\{p^\lambda=p^\lambda_{\le n}+\mathfrak I_{\le n}\}$
and $\{q^\lambda=q^\lambda_{\le n}+\mathfrak I_{\le n}\}$ indexed by compositions $\lambda$ such that $|\lambda |\le n$.
By Propositions~\ref{thm:pmult} and ~\ref{thm:sumsP&QwithM},
the $R$-linear endomorphism $\theta:p^\lambda\to q^\lambda$ is an automorphism.
Applying the automorphism to the equation defining $B_{\lambda,\mu}$ yields
\begin{equation}
\Big(\sum_\eta \bm{\mathfrak b}(\eta)q^\eta\Big) q^\mu=\sum_\lambda B_{\lambda,\mu} q^\lambda.\label{eq:start}
\end{equation}
Now $\sum_\eta \bm{\mathfrak b}(\eta)q^\eta=\sum_\eta \bm{\mathfrak a}(\eta)p^\eta$
by Corollary~\ref{cor:basisdecomp} and $q^\mu=\sum_\beta J_{\beta,\mu} p^\beta$ by Lemma~\ref{lem:inversion}, so
we find
$$
\sum_\gamma \sum_\beta A_{\gamma,\beta}J_{\beta,\mu} p^\gamma
=\sum_\beta \Big(\sum_\eta \bm{\mathfrak a}(\eta)p^\eta\Big) J_{\beta,\mu} p^\beta
=\sum_\alpha \sum_\lambda J_{\alpha,\lambda} B_{\lambda,\mu} p^\alpha.
$$
Comparing the coefficients of $p^\nu$ in both sides gives the result.
\end{proof}

Next we define matrices that will turn out (Theorem~\ref{thm:main}) to be analogues of $M^{(k)}$ for compositions.

\begin{defn}
For each $k>0$, define matrices $\Xi^{(k)}$ and $\Upsilon^{(k)}$
indexed by pairs of compositions of size less than $k$ by
\begin{align*}
\Xi^{(k+1)}&=\begin{pmatrix}\Xi&0\\ 0&\Xi\end{pmatrix} Z^{(k+1)}A^{(k+1)}
\quad\text{and}\quad
\Xi^{(1)}=(1)\\
\Upsilon^{(k+1)}&=Z^{(k+1)}A^{(k+1)} \begin{pmatrix}\Upsilon&0\\ 0&\Upsilon\end{pmatrix}
\quad\text{and}\quad
\Upsilon^{(1)}=(1).
\end{align*}
\end{defn}

The following result is similar to Mathas' conjecture.

\begin{lem}
$\Xi$ and $K\Upsilon K$ are inverse.
\end{lem}

\begin{proof}
Since $B=A^{-1}$ (Lemma~\ref{lem:prodab}) and $Y=Z^{-1}$ (Lemma~\ref{lem:omnibus}),
it follows by induction on $k$ that $\Xi$ is invertible, with inverse
$$
\Xi^{-1}=BY \begin{pmatrix}\Xi^{-1}&0\\ 0&\Xi^{-1}\end{pmatrix}.
$$
We conjugate this equation by $K$.
Lemmas~\ref{matrix1} and~\ref{matrix2} show that $KBYK=KBJ=KJA=ZA$.
Arguing by induction,
$$
\begin{pmatrix}K&-K\\ 0&-K\end{pmatrix} \begin{pmatrix}\Xi^{-1}&0\\ 0&\Xi^{-1}\end{pmatrix} \begin{pmatrix}K&-K\\ 0&-K\end{pmatrix}
=\begin{pmatrix}\Upsilon&0\\ 0&\Upsilon\end{pmatrix}
$$
Hence
$\Upsilon=K\Xi^{-1}K$.
\end{proof}

\section{Increasing elements in products of Jucys-Murphy elements}

\begin{defn}
An element $w\in \mathfrak S_n$ is called {\em increasing} if it has the form
$s_{i_1}s_{i_2}\cdots s_{i_k}$ where $1\le i_1<i_2<\dots <i_k<n$.
To each increasing such $w$, we assign a composition $\phi(w)$ called the \emph{shape} of $w$, as follows.
The set $P=\{s_{i_1},s_{i_2},\dots,s_{i_k}\}$ of simple reflections is totally ordered.
Generate a preorder $\preceq$ by imposing the additional relations $s_{i_j}\sim s_{i_k}$
if they do not commute. This is a composition of $P$.
Let $\phi(w)$ denote the corresponding composition of $k$.
\end{defn}

Conjugacy classes of $\mathfrak S_n$ are commonly indexed by partitions of $n$.
The partition with the same parts as $\phi(w)$ is a partition of $k$, so these are not in general the same.
Let $w$ be increasing of length $k$, and let $\lambda$ be the partition of $k$ with the same parts as $\phi(w)$.
Then the usual shape of $w$ is the partition $(\lambda_1+1,\lambda_2+1,\dots,\lambda_l+1,1,\dots,1)$ of length $n-k$.
In this paper shape means $\phi(w)$ or $\widehat{\phi(w)}$.
For example, $w=s_2 s_3 s_6 s_7s_8\in \mathfrak S_{10}$ has shape $\phi(w)=(2,3)$ rather than $(4,3,1,1,1)$.
There exists an increasing element $w\in \mathfrak S_n$ of shape $\lambda$ iff
$|\lambda |+ \ell(\lambda) \le n$.

In this section, we develop a recurrence (Propositions~\ref{thm:reduction1}
and~\ref{prop:reduction2}) to calculate the bilinear form
$\langle T_w,h\rangle$ for any increasing $w$ and product $h$ of $\LL_1,\dots,\LL_n$.
This is a generalisation of \cite[Prop 2.21]{Mat99} and the proof follows similar lines.

This lemma is inherited from the affine Hecke algebra $\widehat\H$ via the surjective algebra homomorphism $\psi:\widehat\H\to\H$ defined in the Introduction.
\begin{lem}\cite[Lemma 2.15(ii)]{Mat99}\label{lem:mathas}
$$
q\sum_{0\le s<r} \bm{\mathfrak b}(s)\LL_i^s\LL_{i+1}^{r-s}
=T_i\mathcal L_i^rT_i + (q-1)\sum_{1\le s<r} \mathcal L_i^{r-s}T_i\mathcal L_i^{s}.
$$
\end{lem}

\begin{prop}\label{thm:reduction1}
If $w\in \mathfrak S_n$, $h\in \H_n$ and $r$ is a positive integer, then
$$
\langle T_w,h\LL_{n+1}^r\rangle=\sum_{0\le s<r} \langle T_w,ha_s(\LL_1,\dots,\LL_n)\rangle
$$
\end{prop}
\begin{proof}
We define an equivalence relation on $\H_{n+1}$ by
$x\equiv y$ iff $\langle T_w, x\rangle = \langle T_w, y\rangle$ for all $w\in \mathfrak S_n$.
This equivalence relation is preserved by left and right multiplication by $\H_n$.

We prove by induction on $i$ $(0\le i\le n)$ that
\begin{equation}
q^{i-n}T_n\cdots T_{i+1}\LL_{i+1}^rT_{i+1}\cdots T_n
\equiv \sum_{0\le t<r} a_t(\LL_1,\dots,\LL_i)\label{eq:inductive}
\end{equation}
for all positive integers $r$.
The case $i=n$ of \eqref{eq:inductive} is the statement of the theorem,
interpreting the empty product $T_n\cdots T_{i+1}$ as $1$.

First suppose $i=0$.
Then $\LL_{i+1}=1$ so $T_n\cdots T_{i+1}\LL_{i+1}^rT_{i+1}\cdots T_n=q^n\LL_{n+1}\equiv q^n$.
Hence Equation \eqref{eq:inductive} follows.

Let $i>0$ and assume \eqref{eq:inductive} for $i-1$.
By Lemma~\ref{lem:mathas},
\begin{align*}
q\sum_{0\le s<r} \bm{\mathfrak b}(s)\LL_i^sT_n\cdots T_{i+1}\LL_{i+1}^{r-s}T_{i+1}\cdots T_n
=T_n\cdots T_i\LL_i^rT_i\cdots T_n\\ + (q-1)\sum_{1\le s<r} \LL_i^{r-s}T_n\cdots T_{i+1}T_iT_{i+1}\cdots T_n\LL_i^{s}.
\end{align*}
The inductive hypothesis tells us that
$$
T_n\cdots T_i\LL_i^rT_i\cdots T_n
\equiv q^{n+1-i}\sum_{0\le t<r} a_t(\LL_1,\dots,\LL_{i-1})
$$
Furthermore,
$s_n\cdots s_i\cdots s_n=s_i\cdots s_n\cdots s_i\not\in \mathfrak S_n$ so
$$
T_n\cdots T_{i+1}T_iT_{i+1}\cdots T_n\equiv 0.
$$
Therefore (for all positive integers $r$) we have
$$
\sum_{0\le s<r}  \bm{\mathfrak b}(s)\LL_i^sY(r-s)
\equiv \sum_{0\le t<r} a_t(\LL_1,\dots,\LL_{i-1})
$$
where $Y(k)=q^{i-n}T_n\cdots T_{i+1}\LL_{i+1}^kT_{i+1}\cdots T_n$.
Given Lemma~\ref{lem:recurrence}, we may use induction on $k$ to show
$$
Y(k)\equiv \sum_{0\le t<k} a_t(\LL_1,\dots,\LL_i)
$$
This completes the induction on $i$.
\end{proof}

\begin{prop}\label{prop:reduction2}
If $w\in \mathfrak S_n$, $h\in \H_n$ and $r\in \N$, then
$$
\langle T_wT_n,h\LL_{n+1}^r\rangle
= (q-1)\sum_{0\le s<r} \langle T_w,\LL_n^sh\LL_{n+1}^{r-s}\rangle.
$$
\end{prop}
\begin{proof}
If $r=0$, then $\langle T_wT_n,h\rangle=0$ since $ws_n\notin \mathfrak S_n$.
The result follows by induction on $r$ from the calculation:
\begin{align*}
\langle T_wT_n,k\LL_{n+1}\rangle
&=\trace(k\LL_{n+1}T_wT_n)\\
&=\trace(kT_w\LL_{n+1}T_n)\\
&=q^{-1}\trace(kT_wT_n\LL_n T_nT_n)\\
&=\trace(kT_wT_n\LL_n)+(q-1)q^{-1}\trace(kT_wT_n\LL_nT_n)\\
&=\trace(\LL_nkT_wT_n)+(q-1)\trace(kT_w\LL_{n+1})\\
&=\trace(\LL_nkT_wT_n)+(q-1)\trace(k\LL_{n+1}T_w)\\
&=\langle T_wT_n,\LL_nk\rangle + (q-1)\langle T_w,k\LL_{n+1}\rangle
\end{align*}
where $k\in \H_{n+1}$.
\end{proof}

\begin{cor}[Mathas]\label{thm:alt}\cite[Thm 2.7]{Mat99}
If $w\in \mathfrak S_n$ is increasing and $h$ is a product of $k<\ell(w)$ Jucys-Murphy elements,
then $\langle T_w,h\rangle=0$.
\end{cor}

The following is derived as a consequence of Corollary \ref{thm:alt} in \cite{Mat99}.
\begin{cor}[B\"ogeholz~\cite{Bog94}]
Suppose that $1<i_1<i_2<\dots <i_k<n$, and that $w$ is increasing.
Then $\l T_w,L_{i_1}L_{i_2}\dots L_{i_k}\r\neq 0$ if and only if
$w=s_{i_1-1}s_{i_2-1}\dots s_{i_k-1}$.  In this case, $\l T_w,L_{i_1} L_{i_2}\dots L_{i_k}\r=1$.
\label{thm:bogeholz.dec.murphys}
\end{cor}

Propositions~\ref{thm:reduction1} and~\ref{prop:reduction2} provide an (ugly, but) efficient algorithm
for calculating the bilinear form between $T_w$ for $w$ increasing and the quasi-symmetric monomial $\LL^\mu_{\le n}:=p^\mu(\LL_1,\dots,\LL_n)$ for composition $\mu$ with at most $n$ parts.
Technically it is easier if we widen our attention to include incremental terms such as
$$
\LL^{\mu'}_{<n}\LL_{n}^r=\LL^\mu_{\le n}-\LL^\mu_{< n}
$$
where $r=|\mu |-|\mu' |$ is the last part of $\mu$.

\begin{lem}\label{lem:horrid}
Let $w\in \mathfrak S_n$, $\mu$ be a composition of at most $n$ parts and $r$ be a positive integer.
Then we have
\begin{align*}
\langle T_w, \LL^\mu_{\le n}\LL_{n+1}^r\rangle
&=\sum_{ |\lambda |<|\mu | +r}A_{\lambda,\mu} \langle T_w, \LL^\lambda_{\le n}\rangle,\quad\text{and}\\
\langle T_{ws_n}, \LL^\mu_{\le n}\LL_{n+1}^r\rangle&\\
=(q-1)&\sum_{0\le s<r}\sum_{ |\lambda |<|\mu | +r-s}A_{\lambda,\mu}
(\langle T_w, \LL^\lambda_{<n}\LL_n^s\rangle+\langle T_w, \LL^{\lambda'}_{<n}\LL_n^{s+t}\rangle)
\end{align*}
where $t=|\lambda |-|\lambda'|$ is the last part of $\lambda$.
\end{lem}

\begin{proof}
With the notation $\equiv$ from the proof of Proposition~\ref{thm:reduction1},
we have
\begin{align*}
\LL^\mu_{\le n}\LL_{n+1}^r
&\equiv \sum_{0\le s<r} \LL^\mu_{\le n} a_s(\LL_1,\dots,\LL_n)\\
&= \sum_{0\le s<r} \sum_{\substack{|\lambda |=|\mu | + s\\ \ell(\lambda)\le n}} A_{\lambda,\mu} \LL^\lambda_{\le n}\\
&= \sum_{\substack{|\mu |\le |\lambda |<|\mu | + r\\ \ell(\lambda)\le n}} A_{\lambda,\mu} \LL^\lambda_{\le n}\\
&= \sum_{|\lambda |<|\mu | + r} A_{\lambda,\mu} \LL^\lambda_{\le n}
\end{align*}
because $\LL^\lambda_{\le n}=0$ if $\ell(\lambda)>n$ and $A_{\lambda,\mu}=0$ if $|\lambda |<|\mu |$.

Now consider (2). Note that $\langle T_{ws_n}, h\rangle=\langle T_{w}, T_nh\rangle$.
From Proposition~\ref{prop:reduction2} and then part (1), we have
\begin{align*}
T_n\LL^\mu_{\le n}\LL_{n+1}^r
&\equiv (q-1) \sum_{0\le s<r} \LL_n^s \LL^\mu_{\le n}\LL_{n+1}^{r-s}\\
&\equiv (q-1) \sum_{0\le s<r} \LL_n^s \sum_{|\lambda |<|\mu | + r-s} A_{\lambda,\mu} \LL^\lambda_{\le n}.
\end{align*}
Substituting
$$
\LL^\lambda_{\le n}\LL_n^s=\LL^\lambda_{<n}\LL_n^s + \LL^{\lambda'}_{<n}\LL_n^{s+t}
$$
gives the result.
\end{proof}

When $w$ and $\mu$ have the same length,
we can ignore all polynomials in $\LL_i$ of smaller degree.
The recurrence reduces to the one defining $\Xi$.

\begin{thm}\label{thm:main}
If $\lambda$ and $\mu$ are compositions of $k$ and $w\in \mathfrak S_n$ is increasing of shape $\lambda$, then
$$
\l T_w, p^\mu(L_1,\dots,L_n) \r=\Xi_{\lambda',\mu'}^{(k)}.
$$
\end{thm}
\begin{proof}
By Corollary~\ref{thm:alt},
$$
(q-1)^k\l T_w, p^\mu(L_1,\dots,L_n) \r=\l T_w, \LL^\mu_{\le n} \r,
$$
so we shall work with the latter.

In notation of Lemma~\ref{lem:horrid},
the restriction $\ell(w)=|\mu |+r$ makes all terms vanish on the right side.
Therefore (with our notation again) the restriction $\ell(w)=|\mu |$ forces
\begin{equation}
\l T_w, \LL^\mu_{\le n} \r=\l T_{vs_j}, \LL^{\mu'}_{\le j} \LL_{j+1}^r\r \label{firstgo}
\end{equation}
where $w=vs_j$ with $v\in \mathfrak S_j$ and $r=|\mu |-|\mu'|$ is the last part of $\mu$.

We prove by induction on $k$, that for increasing $w$ of length $k$ and compositions $\mu$ such that $|\mu |\le k$,
\begin{equation}
\l T_{ws_n}, \LL^{\mu}_{\le n} \LL_{n+1}^r\r = (q-1)^{k+1} \Xi_{\lambda,\mu}^{(k+1)},\label{inductivehypothesis}
\end{equation}
where $\lambda=\phi(ws_n)'$ and $r = k - |\mu| + 1$.

Assume \eqref{inductivehypothesis} for smaller $k$.
If $w$ is increasing of length $k$ and $\nu$ is a composition such that $|\nu |\le k$,
then we claim that $\l T_w, \LL^\nu_{<n} \LL_n^s\r$ is $(q-1)^k$ times the $(\lambda,\nu)$ coefficient $D_{\lambda,\nu}$ of
$$
D=\begin{pmatrix}\Xi^{(k)}&0\\ 0&\Xi^{(k)}\end{pmatrix}
$$
where $\lambda=\phi(ws_n)'$ and $s=k-|\nu |$.
We have $w\in \mathfrak S_{n-1}$ iff $|\lambda |=k$.
In this case, if $s=0$,
$\l T_w, \LL^\nu_{<n}\r=(q-1)^k\Xi_{\lambda',\nu'}$ by \eqref{firstgo} and \eqref{inductivehypothesis};
if $s>0$, $D_{\lambda,\nu}=0$ by the first recurrence in Lemma~\ref{lem:horrid}.
Alternately suppose $w=vs_{n-1}$ and $|\lambda | < k$.
If $s=0$, then $\l T_w, \LL^\nu_{<n}\r=0$, since the right hand side is in $\mathcal H_{n-1}$
while $w$ lies in the coset $\mathfrak S_{n-1}s_{n-1}$;
if $s>0$,
then $\l T_w, \LL^\nu_{<n}\LL_n^s\r=\l T_{vs_{n-1}}, \LL^\nu_{<n}\LL_n^s\r=(q-1)^k\Xi_{\lambda,\nu}$
by \eqref{inductivehypothesis}.
This proves the claim.

With the notation and hypotheses of \eqref{inductivehypothesis},
consider the second recurrence of Lemma~\ref{lem:horrid}, and ignore all terms which are not
of maximal degree:
\begin{align*}
\l T_{ws_n}, \LL^{\mu}_{\le n}\LL_{n+1}^r \r
&=(q-1)\sum_{ |\nu |\le k}A_{\nu,\mu}
(\langle T_w, \LL^\nu_{<n}\LL_n^s\rangle+\langle T_w, \LL^{\nu'}_{<n}\LL_n^{s+t}\rangle)\\
&=(q-1)^{k+1}\sum_{ |\nu |\le k}A_{\nu,\mu}
(D_{\lambda,\nu}+D_{\lambda,\nu'})\\
&=(q-1)^{k+1}(DZA)_{\lambda,\mu}\\
&=(q-1)^{k+1} \Xi^{(k+1)}_{\lambda,\mu},
\end{align*}
where $s=k-|\nu |$ and $t$ is the last part of $\nu$.
\end{proof}

\section{Conjectures of James and Dipper-James}

For the rest of this paper $R$ denotes an arbitrary commutative ring with $1$,
and $q$ is an invertible element of $R$.
Tensoring with $R$ via the unique ring homomorphism $\Z[q,q^{-1}]\to R:q\mapsto q$,
we find that the equations proven in earlier sections over $\Z[q,q^{-1}]$
are also valid over $R$.

\begin{thm}[James' Conjecture]\label{thm:James.Conj}
Let $k\le n/2$ and consider
the matrix $M^{(k)}$ indexed by partitions $\lambda,\mu\vdash k$ given by
$$
M^{(k)}_{\lambda,\mu}=\l T_w,m_\mu(L_1,\dots,L_n)\r
$$
where $w$ is increasing of shape $\lambda$.
Then $M^{(k)}$ is invertible over $R$.
\end{thm}

\begin{proof}
Fix $k$ and recall the matrix $\Xi$. We need to change the indexation of this matrix.
Let $X$ denote the square matrix indexed by compositions $\lambda, \mu$ of $k$
with entries $X_{\lambda,\mu}=\Xi_{\lambda',\mu'}$ if $k>0$. If $k=0$, $X_{\emptyset,\emptyset}=1$.

Recall from Section 2, that if $\mu$ is a composition, the partition with the same parts as $\mu$ is denoted $\hat\mu$.

Define a matrix $T$ also indexed by compositions $\lambda, \mu$ of $k$ by setting
  \[T_{\lambda,\mu}=\begin{cases}1&\text{if }\lambda=\mu\\ 1&\text{if }\hat\lambda=\mu\\ 0&\text{otherwise.}\end{cases}\]

  Then
  \[\left(XT\right)_{\eta,\mu}=\sum_\lambda X_{\eta,\lambda}T_{\lambda,\mu}=\begin{cases}
\displaystyle\sum_{\hat\lambda=\mu}X_{\eta,\lambda}&\text{if $\mu$ is a partition}\\ X_{\eta,\mu}&\text{otherwise}.\end{cases}
  \]
  Observe that $\sum_{\hat\lambda=\mu}X_{\eta,\lambda}=\l T_w,m_\mu(L_1,\dots,L_n)\r$.

  Now $T$ is invertible, with inverse given by
  \[\left(T^{-1}\right)_{\lambda,\mu}=\begin{cases}1&\text{if }\lambda=\mu\\ -1&\text{if }\hat\lambda
=\mu \text{ but }\lambda\neq\mu\\ 0&\text{otherwise}.\end{cases}\]

  Suppose $\beta$ and $\mu$ are partitions.  Then
  \begin{align*}
    \left(T^{-1}XT\right)_{\beta,\mu}
    &=\sum_\eta \left(T^{-1}\right)_{\beta,\eta}\left(XT\right)_{\eta,\mu}\\
    &=\l T_w,m_\mu(L_1,\dots,L_n)\r
  \end{align*}
  where $w$ is of shape $\eta$.

  Now suppose that $\beta$ is not a partition, but $\mu$ is.  Then
  \begin{align}
    \left(T^{-1}XT\right)_{\beta,\mu}
    &=\sum_\eta \left(T^{-1}\right)_{\beta,\eta}\left(XT\right)_{\eta,\mu}\label{eq:TXT=0}\\
    &=\left(XT\right)_{\beta,\mu}-\left(XT\right)_{\hat\beta,\mu}\notag\\
    &=\l T_{w_\beta},m_\mu(L_1,\dots,L_n)\r-\l T_{w_{\hat\beta}},m_\mu(L_1,\dots,L_n)\r\notag\\
    &=0\notag
  \end{align}
  since $w_\beta$ and $w_{\hat\beta}$ are conjugate and $m_\mu(L_1,\dots,L_n)$ is central.

  The matrix of the Theorem is the submatrix of $T^{-1}XT$ where row and column labels are restricted to partitions.
With that in mind, let $U$ and $V$ be submatrices of $T^{-1}XT$ and $T^{-1}X^{-1}T$ respectively, indexed by partitions.
That is, let $U_{\lambda,\mu}=\left(T^{-1}XT\right)_{\lambda,\mu}$ and
let $V_{\lambda,\mu}=\left(T^{-1}X^{-1}T\right)_{\lambda,\mu}$, for $\lambda$, $\mu$ partitions.
Note that $U$ and $V$ both have entries in $R$.
James' Conjecture states that $U$ is invertible.  We have
  \begin{align*}
    \left(VU\right)_{\eta,\mu}
    &=\sum_{\lambda\text{ a partition}}V_{\eta,\lambda}U_{\lambda,\mu}\\
    &=\sum_{\lambda\text{ a partition}} \left(T^{-1}X^{-1}T\right)_{\eta,\lambda}\left(T^{-1}XT\right)_{\lambda,\mu}\\
    &=\sum_{\lambda} \left(T^{-1}X^{-1}T\right)_{\eta,\lambda}\left(T^{-1}XT\right)_{\lambda,\mu}\\
    \intertext{since $\left(T^{-1}XT\right)_{\lambda,\mu}=0$ if $\lambda$ is not a partition, by \eqref{eq:TXT=0},}
    &=\left(\left(T^{-1}X^{-1}T\right)\left(T^{-1}XT\right)\right)_{\eta,\mu}\\
    &=I_{\eta,\mu}.
  \end{align*}
  Therefore $VU=I$, and hence $U$ is invertible, completing the proof.
\end{proof}

\begin{thm}[The Dipper-James Conjecture]\label{thm:DJconj}
Over a commutative ring $R$ with $1$ and $q\in R$ invertible,
  the set of symmetric functions in Jucys-Murphy elements is the centre of the Hecke algebra $Z(\H)$.
\end{thm}

\begin{proof}
  Consider the matrix $M$ defined for partitions $\lambda$
and $\mu$ by $M^{(k)}_{\lambda,\mu}=\l T_w,m_\mu(L_1,\dots,L_n)\r$ with $w$ of shape $\lambda$.
This is the matrix $U$ of the proof of Theorem~\ref{thm:James.Conj}, but without the restriction that $|\lambda|=|\mu|$,
or that $k\le n/2$.
If $|\lambda|>|\mu|$ then $M_{\lambda,\mu}=0$ by Corollary~\ref{thm:alt}, so $M$ is block triangular,
with rectangular blocks on the diagonal.
Each of these diagonal blocks is $U$ (for a given $|\lambda|=|\mu|$), but with some rows missing.
(A row is missing iff $|\lambda | + \ell(\lambda) > n$.)

  Since $U$ is invertible, each diagonal block has spanning columns.
Therefore $M$ has spanning columns.
It follows that the symmetric polynomials in the Jucys-Murphy elements span the centre of the Hecke algebra.
\end{proof}

We now find a formula for the elementary symmetric functions in Jucys-Murphy elements
in
the
Iwahori-Hecke algebra $\H$ in terms of the Geck-Rouquier basis for the centre,
and hence obtain a corresponding set of generators for the centre of the Hecke algebra,
generalizing a result of Farahat and Higman.

Recall that the {\it $r$'th elementary symmetric function} in $m$ commuting variables
$X_1,\dots,X_m$ is the sum
\[e_r(X_1,\dots,X_m):=\sum_{1\le i_1<i_2<\dots<i_r\le m}X_{i_1}X_{i_2}\dots X_{i_r}.\]

\begin{lem}
If $w$ is increasing then
\[\l T_w,e_r(L_1,\dots,L_n)\r=\begin{cases} 1 &\mbox{if }\ell(w)=r\\
        0   &\mbox{otherwise.}\end{cases}\]
\label{elems}
\end{lem}

\begin{proof}
This is an immediate consequence of Corollary \ref{thm:bogeholz.dec.murphys}.
\end{proof}

Write $\Gamma_\lambda$ for the Geck-Rouquier basis element for the centre of the
Iwahori-Hecke algebra
corresponding to the trace function  $f_{\C_\lambda}$
indexed by partitions of $n$.

\begin{prop} The $r$'th elementary symmetric function in the $n$ Jucys-Murphy elements
is
\[e_r(L_1,\dots,L_n)=\sum_{\ell(\lambda)=n-r}\Gamma_\lambda.\]
\label{thm:elsymfns}
\end{prop}

\begin{proof}
A central element $h$ of $H$ is a linear combination $\sum_\lambda r_\lambda\Gamma_\lambda$ of the Geck-Rouquier basis
with coefficients $r_\lambda$ determined by
$$
\l T_w,h\r =\sum_\lambda r_\lambda \l T_w,\Gamma_\lambda\r =\sum_\lambda r_\lambda f_\lambda(T_w)=r_\mu
$$
where $w$ is an element of minimal length in conjugacy class $\mathcal C_\mu$.
(Note that $\mu$ is not the shape of $w$ in the sense we have used in this paper hitherto.)
Now $h=e_r(L_1,\dots,L_n)$ is central and Lemma~\ref{elems} shows that $r_\lambda=1$ if $\ell(w)=r$, and 0 otherwise.
\end{proof}

We now have an analogue of the following theorem.

\begin{thm}[Farahat-Higman~\cite{FH59}]\label{FH}
The centre of the group algebra $\Z\mathfrak S_n$ is generated as an algebra over $\Z$
by
the set
\[\left\{\sum_{\ell(\lambda)=n-r}\left(\sum_{u\in\C_\lambda}u\right)\Big\vert\ 1\le
r\le
n-1\right\}.\]
\end{thm}

\begin{cor}[to \ref{thm:elsymfns}]
The centre $Z(\H)$ is generated over $R$ by the set
\[\left\{\sum_{\ell(\lambda)=n-r}\Gamma_\lambda\mid 1\le r<n\right\}.\]
\label{cor:gens.for.Z}
\end{cor}

\begin{proof}
  This is immediate from Theorem~\ref{thm:DJconj} and Proposition~\ref{thm:elsymfns}.
\end{proof}

\end{document}